\documentclass{amsart}
\usepackage{amssymb}
\usepackage{amsfonts}
\usepackage{amsmath}
\usepackage[legalpaper,bookmarks=true,colorlinks=true,linkcolor=blue,citecolor=blue]{hyperref}
\usepackage{graphicx}%
\usepackage{fancyhdr}
\usepackage{color}
\usepackage[mathlines]{lineno}
\usepackage{lscape}
\usepackage{epsfig}
\usepackage{natbib}
\usepackage{geometry}
\usepackage{tgbonum}
\usepackage[]{listings}

\fontfamily{qcr}\selectfont

\theoremstyle{plain}

\newtheorem{proposition}{Proposition}

\numberwithin{equation}{section}

\newcommand{\Bin}{\bigskip \noindent}

\newcommand{\Ni}{\noindent}

\begin{document}
\Large
\title{A note on the weak convergence of continuously integrable sequences}
\author{Gane Samb Lo}
\author{Aladji Babacar Niang}

\begin{abstract} A uniformly continuously integrable sequence of real-valued measurable functions, defined on some probability space, is relatively compact in the $\sigma(L^1,L^\infty)$ topology. In this paper, we link such a result to weak convergence theory of bounded measures as exposed in Billingsley (1968) and in Lo(2021) to offer a detailed and new proof using the ideas beneath the proof of prohorov's theorem where the continuous integrability replaces the uniform or asymptotic tightness.\\
 
\noindent $^{\dag}$ Gane Samb Lo.\\
LERSTAD, Gaston Berger University, Saint-Louis, S\'en\'egal (main affiliation).\newline
LSTA, Pierre and Marie Curie University, Paris VI, France.\newline
AUST - African University of Sciences and Technology, Abuja, Nigeria\\
gane-samb.lo@edu.ugb.sn, gslo@aust.edu.ng, ganesamblo@ganesamblo.net\\
Permanent address : 1178 Evanston Dr NW T3P 0J9,Calgary, Alberta, Canada.\\

\noindent $^{\dag \dag}$ Aladji Babacar Niang\\
LERSTAD, Gaston Berger University, Saint-Louis, S\'en\'egal.\\
Email: niang.aladji-babacar@ugb.edu.sn, aladjibacar93@gmail.com\\


\Ni \textbf{AMS 2010 Mathematics Subject Classification:} 60F05; 54C08
\end{abstract}
\maketitle

\section{Introduction}

\Ni The following result is proved, for example in \cite{meyer}.

\begin{proposition} \label{meyer_01}
Let $(f_n)_{n\geq 1}$ be a sequence of measurable mappings:

$$
f_n:\Omega \longrightarrow \mathbb{R},
$$     

\Bin where $\left(\Omega, \mathcal{A}, \mathbb{P}\right)$ is a probability space. Suppose that $(f_n)_{n\geq 0}$ is uniformly continuously integrable (\textit{u.c.i}), that is 

$$
\underset{c\rightarrow +\infty}{\lim}\ \overline{\phi}(c):=\underset{c\rightarrow +\infty}{\lim}\ \underset{n\geq 1}{\sup}\ \int_{\left(\left|f_{n}\right|>c\right)} \left|f_{n}\right| \ d\mathbb{P} = 0.
$$

\Bin Then $(f_n)_{n\geq 1}$ has a sub-sequence $(f_{n_k})_{k\geq 1}$ converging to some $f\in L^1$ in the $\sigma(L^1,L^{\infty})$ topology, i.e.,

$$
\forall g\in L^{\infty}, \ \ \int\ g f_{n_k} \ d\mathbb{P}\rightarrow \int\ g f \ d\mathbb{P}.
$$

\end{proposition}

\Bin Let us make some remarks. \\

\Ni \textbf{Remark 1.} Proposition \ref{meyer_01} implies that: Any  uniformly continuously integrable sequence possesses a subsequence such that the sequence of its mathematical expectations converges in $\mathbb{R}$, since

$$
\mathbb{E}_{\mathbb{P}} (f_{n_k})=\int 1_\Omega \ f_{n_k} \ d\mathbb{P} \rightarrow \int 1_\Omega \ f \ d\mathbb{P}=\mathbb{E}_{\mathbb{P}} (f) 
$$

\Bin This fact is very important in some mathematical problems. For example, it is used in the extension of the Doob decomposition of super-martingales for discrete sequences to time-continuous versions. \\

\Ni \textbf{Remark 2.} Proposition \ref{meyer_01} is very important in some mathematical problems. For example, it used in the extension of the Doob decomposition of super-martingales for discrete sequences to time-continuous versions. \\

\Ni \textbf{Remark 3}. From the weak convergence theory in the sense of \cite{billingsley}, the problem posed here is very similar to prohorov's theorem except we do not have a topology on the domain $\Omega$ of the $f_n$'s. The idea is that we might give a proof following the lines of that of Prohorov's theorem (see Theorem in \cite{ips-wcia-ang}) in which the tightness hypothesis will be replaced by the \textit{u.c.i} condition. The proof given below realizes that idea.

\Ni Let us proceed to the proof.

\section{Proof}

\Ni Let us divide the proof into two parts. \\

\Ni \textbf{Part 1.} Let us suppose that the $f_n$'s are non-negative. Let us set measures $\nu_n$:

$$
d\nu_n = f_n d\mathbb{P}, \ n\geq 1.
$$ 

\Bin Let 

$$
\mathcal{C}_0 = \left\{g\in L^{\infty}, \ \left\|g\right\|\leq 1\right\}.
$$

\Bin Since the $f_n$'s are \textit{u.c.i}, we have for $n\geq 1$, for $c>0$,

\begin{eqnarray*}
\int |f_n| \ d\mathbb{P} &=& \int_{(|f_n|\leq c)} |f_n| \ d\mathbb{P} + \int_{(|f_n|>c)} |f_n| \ d\mathbb{P} \\
&\leq& c + \underset{n\geq 1}{\sup}\ \int_{\left(\left|f_{n}\right|>c\right)} \left|f_{n}\right| \ d\mathbb{P} \\
&=& c + \overline{\phi}(c) 
\end{eqnarray*}

\Ni and 

$$
\overline{\phi}(c)\rightarrow 0, \ as \ c\rightarrow +\infty.
$$

\Bin So, for some $c_0>0$, we have

$$
\underset{n\geq 1}{\sup}\ \int |f_n| \ d\mathbb{P} \leq c_0 + 1 =: C.
$$

\Bin Now, for $g\in \mathcal{C}_0$, $\forall n\geq 1$,

$$
\left|\int g f_n \ d\mathbb{P}\right| \leq \int |f_n| \ d\mathbb{P} \leq C.
$$

\Bin Hence 

$$
\overline{\nu}_n = \left(\nu_n(g), \ g\in \mathcal{C}_0\right)\in \left[0, \ C\right]^{\mathcal{C}_0}.
$$

\Bin Since $\left[0, \ C\right]^{\mathcal{C}_0}$ as an arbitrary product of the compact set $\left[0, \ C\right]$, by Tchychonov's theorem, we may extract  a sub-sequence $(\overline{\nu}_{n_k})_{k\geq 1}$ from $(\overline{\nu}_n)_{n\geq 1}$ converging to $\overline{\nu} = \left(\nu(g), \ g\in \mathcal{C}_0\right)$, that is

$$
\forall g\in \mathcal{C}_0, \ \nu_{n_k}(g)\rightarrow \nu(g), 
$$

\Bin with 

$$
\forall g\in \mathcal{C}_0, \ \nu(g)\leq C. 
$$

\Bin We have to show that $\nu$ is still a measure on $\left(\Omega, \mathcal{A}, \mathbb{P}\right)$ for $g=1_A$, \ $A\in \mathcal{A}$. But in general, we have for $(\alpha,\beta)\in \mathbb{R}^2$: $|\alpha|+|\beta|\leq 1$, \ $(g,h)\in \mathcal{C}_0^2$,

\begin{eqnarray*}
\nu_{n_k}(\alpha g + \beta h) &=& \alpha \nu_{n_k}(g) + \beta \nu_{n_k}(h) \\
&\rightarrow& \alpha \nu(g) + \beta \nu(h) \ as \ k\rightarrow +\infty
\end{eqnarray*}

\Ni and since $\alpha g + \beta h\in \mathcal{C}_0$, we have

$$
\nu_{n_k}(\alpha g + \beta h)\rightarrow \nu(\alpha g + \beta h) \ as k\rightarrow +\infty
$$

\Bin and hence 

$$
\nu(\alpha g + \beta h) = \alpha \nu(g) + \beta \nu(h).
$$

\Bin We prove similarly that for $(g_1,\cdots,g_p)\in \mathcal{C}_0^p$, \ $p\geq 2$, for $(\alpha_1,\cdots,\alpha_p)\in \mathbb{R}^p$: $|\alpha_1| +\cdots+|\alpha_p|\leq 1$,

$$
\nu_{n_k}(\alpha_1 g_1 +\cdots+ \alpha_p g_p)\rightarrow \alpha_1 \nu(g_1) +\cdots+ \alpha_p \nu(g_p).
$$

\Bin Now for any pair-wise disjoint measurable sets $(A_j)_{1\leq j\leq p}$, $1_{\sum_{j=1}^{p} A_j}\in \mathcal{C}_0$, and

$$
\nu_{n_k}\left(\sum_{j=1}^{p} A_j\right)=:\nu_{n_k}\left(1_{\sum_{j=1}^{p} A_j}\right)\rightarrow \nu\left(1_{\sum_{j=1}^{p} A_j}\right) 
$$

\Bin and 

$$
\nu_{n_k}\left(\sum_{j=1}^{p} A_j\right)=\sum_{j=1}^{p} \nu_{n_k} (1_{A_j})\rightarrow \sum_{j=1}^{p} \nu(1_{A_j}).
$$

\Bin So,

$$
\nu\left(\sum_{j=1}^{p} A_j\right) = \sum_{j=1}^{p} \nu(A_j).
$$

\Bin Let us prove that $\nu$ is a measure by extending the finite additivity to a $\sigma$-additivity. \\

\Ni But since $\nu$ is non-negative, finite and additive, so it is enough to show that $\nu$ is continuous at $\emptyset$:

$$
\forall A_p\downarrow \emptyset, \ \nu(A_P)\downarrow 0 \ as \ p\uparrow +\infty.
$$ 

\Bin By Theorem \cite{ips-mfpt-ang}, an \textit{u.c.i} sequence is also uniformly absolutely:

$$
\forall \varepsilon>0, \ \biggr(\exists \delta>0, \ \forall A\in \mathcal{A}: \mathbb{P}(A)< \delta, \ \underset{n\geq 1}{\sup}\ \int_{A} |f_n| \ d\mathbb{P}<\varepsilon\biggr).
$$

\Bin Let $A_p\downarrow \emptyset$, \ $A_p\in \mathcal{A}$. So there exists $p_0\in \mathbb{N}$, such that for any $p\geq p_0$, $\mathbb{P}(A_p)<\epsilon$ and hence

$$
\underset{k\geq 1}{\sup}\ \int_{A_p} |f_{n_k}| \ d\mathbb{P}<\epsilon,
$$

\Bin that is,

$$
\forall p\geq p_0, \ \forall k\geq 1, \ |\nu_{n_k}(A_p)|\leq \int_{A_p} |f_{n_k}| \ d\mathbb{P}<\varepsilon.
$$

\Bin By letting $k\rightarrow +\infty$, we get

$$
\forall p\geq p_0, \ \nu(A_p)<\epsilon.
$$

\Bin We proved that 

$$
\nu(A_P) \rightarrow 0 \ as \ p\rightarrow +\infty.
$$ 

\Bin $\nu$ is non-decreasing as an additive and non-negative mapping on $\mathcal{A}$. So $\nu$ is a measure on $\left(\Omega, \mathcal{A}, \mathbb{P}\right)$, Let us show that $\nu$ is continuous with respect to $\mathbb{P}$, that is

$$
d\nu = f \ d\mathbb{P}, \ f\in L^1.
$$

\Bin We have for $A\in \mathcal{A}$ such that $\mathbb{P}(A)=0$,

$$
0=\int_{A} f_{n_k} \ d\mathbb{P} = \nu_{n_k}(A) \rightarrow \nu(A)
$$

\Bin and then $\nu(A)=0$. So $\nu$ is continuous with respect to $\mathbb{P}$. By Random Nikodym theorem, there exists $f\in L^1$ (since $\nu$ is finite) such that

$$
d\nu = f \ d\mathbb{P}.
$$ 

\Bin So 

\begin{equation}
\forall A\in \mathcal{A}, \ \int 1_{A} \ f_{n_k} \ d\mathbb{P}\rightarrow \int 1_{A} \ f \ d\mathbb{P}. \label{wc_01}
\end{equation}

\Bin Now, \eqref{wc_01} is automatically holds if $1_A$ is replaced by an elementary function on $(\Omega, \ \mathcal{A}, \mathbb{P})$. Next any element $g$ of $\mathcal{C}_0$ is limit of a sequence $(g_p)$ of elementary functions with $|g_p|\leq |g|\leq 1$. So for $p\geq 1$, $k\geq 1$, $\varepsilon>0$, we have  

\begin{eqnarray*}
&&\left|\int g \ f_{n_k} \ d\mathbb{P} -\int g \ f \ d\mathbb{P}\right| \label{fxx}\\
&&\leq \left|\int (g-g_p) \ f_{n_k} \ d\mathbb{P}\right| +\left|\int g_p \ f_{n_k} \ d\mathbb{P} -\int g \ f \ d\mathbb{P}\right|\\
&&\leq \left|\int_{(|g-g_p|\leq \varepsilon)} (g-g_p) \ f_{n_k} \ d\mathbb{P}\right|+\left|\int_{(|g-g_p|>\varepsilon)} (g-g_p) \ f_{n_k} \ d\mathbb{P}\right| +\left|\int g_p \ f_{n_k} \ d\mathbb{P} -\int g \ f \ d\mathbb{P}\right| \notag\\
&&\leq \varepsilon \int \ f_{n_k} \ d\mathbb{P}+\int 1_{(|g-g_p|>\varepsilon)} \ f_{n_k} \ d\mathbb{P}
+\left|\int g_p \ f_{n_k} \ d\mathbb{P} -\int g \ f \ d\mathbb{P}\right|. \notag
\end{eqnarray*}

\Bin By letting $p\geq 1$ and $\varepsilon>0$ fixed, and letting $k\rightarrow +\infty$, we get

\begin{eqnarray*}
\limsup_{k\rightarrow +\infty}\left|\int g \ f_{n_k} \ d\mathbb{P} -\int g \ f \ d\mathbb{P}\right|
\leq \varepsilon \int \ f \ d\mathbb{P}+\int 1_{(|g-g_p|>\varepsilon)} \ f \ d\mathbb{P}
+\left|\int g_p \ f \ d\mathbb{P} -\int g \ f \ d\mathbb{P}\right|. \notag
\end{eqnarray*}

\Bin Now, we let $p\rightarrow +\infty$, to get

$$
\left|\int g_p \ f \ d\mathbb{P} -\int g \ f \ d\mathbb{P}\right| \rightarrow 0,
$$

\Bin by the dominated convergence theorem applied to $L^1 \ni (2f) \geq |g-g_p|\rightarrow 0$ and 

$$
\int 1_{(|g-g_p|>\varepsilon)} \ f \ d\mathbb{P} \rightarrow 0,
$$

\Bin since $g_p \rightarrow g$ and hence $g_p \rightarrow_{\mathbb{P}} g$ and next $\mathbb{P}(|g-g_p|>\varepsilon)\rightarrow 0$ and finally the integrable function $f$ is continuous. We get for any $\varepsilon>0$

\begin{eqnarray*}
\limsup_{k\rightarrow +\infty}\left|\int g \ f_{n_k} \ d\mathbb{P} -\int g \ f \ d\mathbb{P}\right|\leq \varepsilon \int f \ d\mathbb{P}.
\end{eqnarray*}

\Bin We conclude by letting $\varepsilon\downarrow 0$ and get

$$
\forall g \in \mathcal{C}_0, \int g \ f_{n_k} \ \ d\mathbb{P} \rightarrow \int g \ f \ \ d\mathbb{P}.
$$

\Bin Finally, for any $g\in L^{\infty}$, $\left\|g\right\|_{\infty}\geq 1$,

$$
\overline{g}=\frac{g}{\left\|g\right\|_{\infty}}\in \mathcal{C}_0
$$

\Bin and

$$
\int g f_{n_k} \ d\mathbb{P} = \left\|g\right\|_{\infty} \int \overline{g} f_{n_k} \ d\mathbb{P}\rightarrow \left\|g\right\|_{\infty} \int \overline{g} f \ d\mathbb{P} = \int g f \ d\mathbb{P}.
$$

\Bin The proof of \textbf{Part 1} is complete. \\

\Ni \textbf{Part 2}. Here the $f_n$s are not necessarily non-negative.  We define for $g\in \mathcal{C}_0$

\begin{eqnarray*}
\nu_n(g) &=& \int g f_n^{+} d\mathbb{P} - \int g f_n^{-} \ d\mathbb{P} \\
&=& \nu_n^{+}(g) - \nu_n^{-}(g).
\end{eqnarray*}

\Bin The sequences $(f_n^{+})_{n\geq 1}$ and $(f_n^{-})_{n\geq 1}$ are still \textit{u.c.i} since for example

$$
\int_{(f_n^{+}>c)}  f_n^{+} \ d\mathbb{P}\leq \int_{(|f_n|>c)}  f_n^{+} \ d\mathbb{P}\leq \int_{(|f_n|>c)}  |f_n| \ d\mathbb{P}\leq \overline{\phi}(c).
$$

\Bin So for some constant $C>0$, we have 

$$
0\leq \nu_n^{+}\leq C \ and \ 0\leq \nu_n^{-}\leq C.
$$

\Bin By \textbf{Part 1}, we extract a subsequence $(\nu_{n^\prime}^{+})_{n^\prime\geq 1}$ from $(\nu_{n}^{+})_{n\geq 1}$ that weakly converges to
$\nu^+$ with $d\nu^+=f^{+} d\mathbb{P}$, $f^{+}\in L^1$. Next From  $(\nu_{n^\prime}^{-})_{n^\prime\geq 1}$ we extract a subsequence $(\nu_{n_k}^{-})_{k\geq 1}$ and $(\nu_{n_k}^{-})_{k\geq 1}$ that weakly converges to $\nu^-$ with $d\nu^-=f^{-} d\mathbb{P}$, $f^{-}\in L^1$. Hence

$$
\nu_{n_k}^{+} \rightarrow_{\sigma(L^1,L^\infty)} \nu^{+} \ \ and \ \ \nu_{n_k}^{-} \rightarrow_{\sigma(L^1,L^\infty)} \nu^-, 
$$

\Bin and

$$
\nu_{n_k}=\nu_{n_k}^{+}-\nu_{n_k}^{-} \rightarrow_{\sigma(L^1,L^\infty)} \nu^{+}-\nu^{-}=\nu, 
$$

\Bin with $d\nu=f \ d\mathbb{P}$, with $f=f^{+}-f^{-} \in L^1$. $\blacksquare$

\end{document}